\documentclass[11pt]{amsart}          
\usepackage{amsfonts}
\usepackage{amssymb}  
\usepackage{amsthm} 
\usepackage{amsmath} 
\usepackage{tikz-cd}
\usepackage{float}
\usepackage{scalerel}
\usepackage[]{hyperref}
\hypersetup{
  colorlinks,
  linkcolor=blue,
  citecolor=blue,
  urlcolor=blue}
\newcommand{\cat}[1]{\mathcal{#1}}
\newcommand{\Cat}[1]{\mathbf{#1}}
\newcommand{\Set}{\Cat{Set}}

\author{Bartosz Milewski}
\title{Compound Optics}
\begin{document}
\begin{abstract}
Simple optics are defined using actions of monoidal categories. Compound optics arise, for instance, as natural transformations between polynomial functors. Since a monoidal category is a special case of a bicategory, we formulate complex optics using the action of a bicategory. We show that polynomial optics are a special case of complex optics defined by the action of bicategory $\mathbf{Prof}$ on co-presheaves. 

\end{abstract}
\maketitle{}
\date{}

\section{Introduction}

Optics, as it is used in computer science, provides a way of focusing on and modifying individual parts of data structures. In category theory, it provides a window into the inner structure of objects. In general, optics tries to decompose an object into a focus, or multiple foci, and a transformation that relates the foci to the whole object. In simple optics the foci are objects in the same category as the whole object. This framework was first described in \cite{milewski17}, formalized in \cite{riley18}, and extended in \cite{clarke22}.

Recently, new types of optics arose in the study of polynomial functors \cite{milewski21a, milewski21b}. These new compound optics require a generalization of the standard framework by grouping the foci into separate categories. The contribution of this paper is to show that transformations that relate these foci to the composite objects, rather than forming an actegory, require the use of a bicategory. Composition of compound optics is explained in terms of Kan extensions \cite{maclane78}. 

\section{Simple Optics}

Simple optics is built from several layers. At the bottom we have actegories that are defined by the action of a monoidal category $\cat M$ upon regular categories. Such an action is often written using infix notation $m \bullet c$, where:
\[ \bullet \colon \cat M \times \cat C \to \cat C \]
or, in the curried form reminiscent of the Cayley's theorem, as a monoidal functor:
\[ \bullet \colon \cat M \to [\cat C, \cat C] \]
Two such actions in two categories:
\begin{align*}
\bullet_{\scaleto{1}{4pt}} &\colon \cat M \to [\cat C, \cat C ] \\
 \bullet_{\scaleto{2}{4pt}} &\colon \cat M \to [\cat D, \cat D ] 
\end{align*}
can be combined into a mixed-variance action:
\[ \bullet \colon \cat M^{op} \times \cat M \to [\cat C^{op} \times \cat D, \cat C^{op} \times \cat D] \]
where:
\[  \langle m, n  \rangle \bullet \langle c, d \rangle = \langle m \bullet_{\scaleto{1}{4pt}} c, n \bullet_{\scaleto{2}{4pt}} d \rangle \]
The hom-sets in $\cat C^{op} \times \cat D$ are then extended (later we'll see that these are in fact Kan extensions) along these actions:
\[ (\cat C^{op} \times \cat D)(  \langle m, n  \rangle \bullet \langle a, b \rangle, \langle s, t \rangle) \]
Finally, after averaging (by taking the coend, \cite{loregian19}) over the monoidal category $\cat M$, we get the most general simple (mixed) optic:
\[ \mathcal{O}\langle a, b\rangle \langle s, t \rangle = \int^{m \colon \cat M} (\cat C^{op} \times \cat D)(  \langle m, m  \rangle \bullet \langle a, b \rangle, \langle s, t \rangle) \]
Or, in the more familiar form \cite{clarke22}:
\[ \mathcal{O}\langle a, b\rangle \langle s, t \rangle = \int^{m \colon \cat M} \cat C (s, m \bullet_{\scaleto{1}{4pt}} a) \times \cat D ( m \bullet_{\scaleto{2}{4pt}} b, t) \]
Here, $a$ and $b$ are interpreted as the foci; and $s$ and $t$ are the whole objects.

These optics are then interpreted as hom-sets in the category $\Cat {Opt}$ \cite{riley18} where objects are pairs like $\langle a, b \rangle$ and $\langle s, t \rangle$. 

Not all of optics fall into the category of simple optics. For one, it's possible to compose disparate optics, like lenses with prisms to get new optics. Then there are optics that arise as natural transformations between polynomial functors. These optics have a representation as ``ommatidia'' \cite{milewski21b}, which is similar to simple optics, but goes beyond the representation of monoidal categories. A larger framework is needed.

\section{Bicategories}

A bicategory $\mathbb{C}$ has objects as 0-cells, arrows as 1-cells, and arrows between arrows as 2-cells. This allows for the relaxation of the usual commuting rules for diagrams. For instance, composing a 1-cell with an identity 1-cell  may result in a 1-cell that's not exactly equal to the original but is related to it by an invertible 2-cell called the unitor. Similarly, associativity of 1-cell composition may be defined up to a 2-cell isomorphism called the associator. These have to satisfy additional coherence conditions.

It's worth pointing out that 1-cells between any pair of 0-cells $i$ and $j$ form a hom-category $\mathbb C (i, j)$.

If the axioms of a category are satisfied on the nose, a bicategory is called a 2-category. The common example of a 2-category is the category of small categories $\Cat{Cat}$, with categories as 0-cells, functors as 1-cells, and natural transformations as 2-cells. 1-cells between any two categories $\cat C$ and $\cat D$ form a functor category 
\[ \Cat{Cat} (\cat C, \cat D) =  [\cat C, \cat D] \]

An important example of a bicategory is a single-object bicategory which is equivalent to a monoidal category. Its 1-endo-cells represent the objects of the monoidal category: their composition is the monoidal tensor product. The identity 1-cell is the unit object $I$. The 2-cells represent the morphism of the monoidal category. 

A pseudofunctor $F$ between bicategories $\mathbb C$ and $\mathbb D$ is a mapping that weakly preserves the bicategorical structure. In particular its action on 1-cells is a regular functor between hom-categories:
\[ \mathbb C (i, j) \to \mathbb D (F i, F j) \]
A pseudofunctor preserves composition and identity up to 2-cell isomorphism. There are also some additional coherence condition.

A pseudofunctor $F$ from a single-object bicategory $\Cat M$ to $\Cat{Cat}$ defines an actegory, that is an action of a monoidal category on some category $\cat C$. Indeed, such a functor maps the 0-cell $*$ to a category $\cat C$, and it maps (functorially) the hom-category $\Cat M(*, *)$ to the endo-functor category $[\cat C, \cat C]$.

Recall that the 1-cells in $\Cat M(*, *)$ are objects of the monoidal category and their composition is the tensor product. Thus if we have 1-cells $m$ and $n$ in $\Cat M(*, *)$, $F$ maps their product $m \otimes n$ to the composition of endofunctors $F m \circ F n$. Pseudo-functoriality means that unit and associativity are weakly preserved.

\section{Simple Optic Explained Using Bicategories}

A simple optic can be defined using a one-object bicategory $\Cat M$ and a pair of pseudofunctors into $\Cat{Cat}$. These functors pick a pair of categories $\cat C$ and $\cat D$ and a pair of actions:
\begin{align*}
\bullet_{\scaleto{1}{4pt}} &\colon \Cat M(*, *) \to [\cat C, \cat C ] \\
 \bullet_{\scaleto{2}{4pt}} &\colon \Cat M(*, *) \to [\cat D, \cat D ] 
\end{align*}
These two actions can be combined into one:
\[ \bullet \colon \Cat M(*, *)^{op} \times \Cat M(*, *) \to [\cat C^{op} \times \cat D, \cat C^{op} \times \cat D] \]
where:
\[  \langle m, n  \rangle \bullet \langle c, d \rangle = \langle m \bullet_{\scaleto{1}{4pt}} c, n \bullet_{\scaleto{2}{4pt}} d \rangle \]

\section{Kan Extensions}

To explain how the category of optics arises from monoidal actions, it's instructive to go through an intermediate construction that involves Kan extensions. A point-wise left Kan extension of a co-presheaf $F$ along a functor $P$ can be written as a coend:
\[ (\text{Lan}_P F) b \cong \int^{c} \cat B (P c, b) \times F c \]
The categories and functors in this definition are shown in the following diagram:
\[
 \begin{tikzcd} [row sep=1cm, column sep=1cm]
 \cat C
 \arrow[rr, "F", "" {name=ID, below} ]
 \arrow[d, bend right, "P"']
 && \Set
 \\
 \cat B
  \arrow[urr, bend right, "\text{Lan}_P\;F"']
 \end{tikzcd}
\]

Kan extensions are composable. An extension along a composition of functors is a composition of Kan extensions:
\[ \text{Lan}_{Q\circ P} F \cong  \text{Lan}_Q (\text{Lan}_P F) \]
Pictorially:
\[
 \begin{tikzcd} [row sep=1cm, column sep=1cm]
 \cat C
 \arrow[rr, "F", "" {name=ID, below} ]
 \arrow[d, bend right, "P"']
 && \Set
 \\
 \cat B
  \arrow[urr, bend right, "\text{Lan}_P F"]
 \arrow[d, bend right, "Q"']
 \\
 \cat B'
  \arrow[uurr, bend right, "\text{Lan}_Q\; (\text{Lan}_P\; F)"']
 \end{tikzcd}
\]

Representable functors form a basis of the co-presheaf category (any co-presheaf is a colimit of representables). It turns out that Kan extensions of representables are also composable. 

\section{Composition on Representables}

Kan extensions of representables can be used to define profunctors that are parameterized by functors, which can be seen as generalized hom-sets between categories, gluing them into collages. Consider a functor $P \colon \cat C \to \cat D$ and two objects $c \colon \cat C$ and $d \colon \cat D$. Define the profunctor:
\[ \Pi_P (c, d) =  \left(\text{Lan}_P \mathcal{Y}_c \right) d = \int^{c'} \cat D(P c', d) \times \cat C(c, c') \]
where 
\[ (\mathcal{Y}_c) c' = \cat C (c, c') \]
is the Yoneda functor.
\[
 \begin{tikzcd} [row sep=1cm, column sep=1cm]
 \cat C
 \arrow[rr, "\mathcal Y_c", "" {name=ID, below} ]
 \arrow[d, bend right, "P"']
 && \Set
 \\
 \cat B
  \arrow[urr, bend right, "{\Pi_P(c, -)}"']
 \end{tikzcd}
\]

For any two composable functors $P$ and $Q$ we can define the composition of the profunctors:
\[ \Pi_Q (c, d) \times \Pi_P (d, e) \to \Pi_{Q \circ P} (c, e) \]
Using the definition of Kan extension, together with the fact that products in $\Set$ preserve coends, and the Fubini theorem, we get:
\[ \left( \text{Lan}_P \mathcal{Y}_c \right) d \times \left( \text{Lan}_Q \mathcal{Y}_d \right) e \cong 
    \int^{c', d'}  \cat D (P c', d) \times \cat C (c, c') \times \cat E (Q d', e) \times \cat D (d, d') \]
 Two of the hom-sets are composable:
\[ \circ \colon \cat D (d, d') \times \cat D (P c', d)  \to D (P c', d') \]
which gives us the desired mapping to:
\[ \int^{c', d'}  \cat D (P c', d') \times \cat C (c, c') \times \cat E (Q d', e) \times \cong
   \left( \text{Lan}_Q\; (\text{Lan}_P\; \cat Y_c) \right) e \cong
  \left( \text{Lan}_{Q \circ P} \cat Y_c \right) e \]

Later we'll use this formula with the functor $P$ replaced by the action $m \bullet$. 

\section{Averaging over Actions}

In this picture, simple optics are recovered by considering co-presheaves on the product category $\cat C^{op} \times \cat D$. We construct the left Kan extension of the representable:
\[ \mathcal Y_{\langle a, b \rangle} \langle c, d \rangle = (\cat C^{op} \times \cat D)(\langle a, b \rangle, \langle c, d \rangle) = \cat C(c, a) \times \cat D(b, d) \]
along the double action of the monoidal category:
\[ \bullet \colon \cat M^{op} \times \cat M \to [\cat C^{op} \times \cat D, \cat C^{op} \times \cat D] \]
as illustrated in the following diagram:
\[
 \begin{tikzcd} [row sep=1cm, column sep=1cm]
 \cat C^{op} \times \cat D
 \arrow[rr, "\mathcal Y_{\langle a, b \rangle}", "" {name=ID, below} ]
 \arrow[d, bend right, "{\langle m, n \rangle \bullet}"']
 && \Set
 \\
 \cat{C}^{op} \times \cat D
  \arrow[urr, bend right, "{\text{Lan}_{\langle m, n \rangle \bullet}\; \mathcal Y_{\langle a, b \rangle}}"']
 \end{tikzcd}
\]
This gives us the integrand of the simple optic:
\[ \left(\text{Lan}_{\langle m, n \rangle \bullet}\; \mathcal Y_{\langle a, b \rangle}\right) \langle s, t \rangle \cong \int^{\langle c, d \rangle \colon \cat C^{op} \times \cat D} \cat (\cat C^{op} \times \cat D)(\langle m, n \rangle \bullet \langle c, d \rangle, \langle s, t \rangle) \times \mathcal Y_{\langle a, b \rangle} \langle c, d \rangle \]
The coend can be reduced using the co-Yoneda lemma, resulting in:
\[ (\cat C^{op} \times \cat D)(\langle m, n \rangle \bullet \langle a, b \rangle, \langle s, t \rangle)\]
which can be written in the more familiar form as:
\[ \cat C (s, m \bullet_{\scaleto{1}{4pt}} a) \times \cat D (n \bullet_{\scaleto{2}{4pt}} b, t) \]

The optics themselves are obtained by averaging over the actions, that is taking the coend over the monoidal category $\cat M$:
\[  \mathcal{O}\langle a, b\rangle \langle s, t \rangle \cong \int^{m \colon \cat M} \left(\text{Lan}_{\langle m, m \rangle \bullet}\; \mathcal Y_{\langle a, b \rangle}\right) \langle s, t \rangle \]

To see that optics are composable, let's rewrite them using the composable profunctors $\Pi$:
\[  \mathcal{O}\langle a, b\rangle \langle s, t \rangle \cong \int^{m \colon \cat M} \Pi_{\langle m, m \rangle \bullet}\langle a, b\rangle \langle s, t \rangle\]
We have to construct a mapping:
\[ \int^{m, n} \Pi_{\langle m, m \rangle \bullet}\langle a, b\rangle \langle s, t \rangle \times \Pi_{\langle n, n \rangle \bullet}\langle s, t\rangle \langle s', t' \rangle 
\to \int^{k} \Pi_{\langle k, k \rangle \bullet}\langle a, b\rangle  \langle s', t' \rangle\]
We first compose the two $\Pi$'s to get:
\[ \int^{m, n} \Pi_{\langle m \otimes n, m \otimes n \rangle \bullet}\langle a, b\rangle  \langle s', t' \rangle\]
where the tensor product is the result of the composition of two actions $n \circ m$.

The remaining step is to use the universal property of the coend with the injections:
\[ \iota_{\langle m \otimes n, m \otimes n \rangle} \colon \Pi_{\langle m \otimes n, m \otimes n \rangle \bullet}\langle a, b\rangle  \langle s', t' \rangle \to \int^{r} \Pi_{\langle r, r \rangle \bullet}\langle a, b\rangle  \langle s', t' \rangle \]

\section{Compound Optics}
To describe compound optics we generalize the monoidal category of simple optics to a bicategory $\mathbb{M}$. The monoidal action is then replaced by a pseudofunctor $F \colon \mathbb M \to \Cat{Cat}$. Such a pseudofunctor maps 0-cells to categories. A hom-category between two 0-cells is mapped to a functor category:
\[ \mathbb{M} (i, j) \to [\cat C, \cat C'] \]
where $\cat C = F i$ and $\cat C' = F j$. 
Functoriality means that the composition of 1-cells is mapped to composition of functors. Notice that this time we are not talking about endofunctors. 

To define compound optics we take two such pseudofunctors that define a pair of actions:
\begin{align*}
\bullet_{\scaleto{1}{4pt}} &\colon \mathbb M(i, j) \to [\cat C, \cat C' ] \\
 \bullet_{\scaleto{2}{4pt}} &\colon \mathbb M(i, j) \to [\cat D, \cat D' ] 
\end{align*}
These two actions can be combined into a mixed-variance action:
\[ \bullet \colon \mathbb M(i, j)^{op} \times \mathbb M(i, j) \to [\cat C^{op} \times \cat D, \cat {C'}^{op} \times \cat {D'}] \]
If $m$ and $m'$ are elements of $\mathbb M(i, j)$, their action is defined as:
\[  \langle m, m'  \rangle \bullet \langle c, d \rangle = \langle m \bullet_{\scaleto{1}{4pt}} c, m' \bullet_{\scaleto{2}{4pt}} d \rangle \]
This time the result of the action lives in a different category than the source. 

We use the left Kan extension to extend the representable:

\[ \mathcal Y_{\langle a, b \rangle} \langle c, d \rangle = (\cat C^{op} \times \cat D)(\langle a, b \rangle, \langle c, d \rangle) = \cat C(c, a) \times \cat D(b, d) \]
along the action of $\langle m, m' \rangle$:
\[
\begin{tikzcd} [row sep=1cm, column sep=1cm]
i 
\arrow[d, "{m, m'}"]
&& \cat C^{op} \times \cat D
 \arrow[rr, "\mathcal Y_{\langle a, b \rangle}", "" {name=ID, below} ]
 \arrow[d, bend right, "{\langle m, m' \rangle \bullet}"']
 && \Set
 \\
 j && \cat{C'}^{op} \times \cat D'
  \arrow[urr, bend right, "{\text{Lan}_{\langle m, m' \rangle \bullet}\; \mathcal Y_{\langle a, b \rangle}}"']
 \end{tikzcd}
\]
We've seen earlier that these Kan extensions are composable as long as the actions are composable (using regular functor composition). This time the profunctors:
\[ \Pi_{\langle m, m' \rangle \bullet}\langle a, b\rangle \langle s, t \rangle = \left( \text{Lan}_{\langle m, m' \rangle \bullet}\; \mathcal Y_{\langle a, b \rangle} \right) \langle s, t \rangle \]
connect objects in four different categories, not counting the hom-category $\mathbb M (i, j)$.

The final step is averaging over $\mathbb M (i, j)$. We define the composite optic as a coend:

\[  \mathcal{O}_{i, j} \langle a, b\rangle \langle s, t \rangle \cong \int^{m \colon \mathbb M (i, j)} \left(\text{Lan}_{\langle m, m \rangle \bullet}\; \mathcal Y_{\langle a, b \rangle}\right) \langle s, t \rangle \]

These optics compose as long as the corresponding 1-cells in $\mathbb M$ compose. Consider two composable 1-cells, $n \colon \mathbb M (i, j)$ and $m \colon \mathbb M (j, k)$. Their composition is $m \circ n \colon \mathbb M (i, k)$. The composition of compound optics follows the same pattern as that of simple optics, except that the composition of actions no longer reduces to a tensor product. The injection that completes the proof is:

\[ \iota_{\langle m \circ n, m \circ n \rangle} \colon \Pi_{\langle m \circ n, m \circ n \rangle \bullet}\langle a, b\rangle  \langle s', t' \rangle 
\to \int^{r \colon \mathbb M(i, k)} \Pi_{\langle r, r \rangle \bullet}\langle a, b\rangle  \langle s', t' \rangle \]
The identity optic is a member of:
\[ \mathcal{O}_{i, i} \langle a, b\rangle \langle a, b \rangle \]
generated by the identity 1-cell in $\mathbb M (i, i)$.

\section{Polynomial optics}

A practical example of compound optics arises in the study of polynomial functors \cite{gambino_kock_2013}. Such functors can be written as coproducts of representables:
\[ P(y) = \sum_{n \in N} s_n \times \mathbf{Set}(t_n, y) \]
where $s_n$ and $t_n$ are families of sets indexed by the elements of the set $N$. A natural transformation from $P$ to $Q$, where $Q$ is given by the formula:
\[ Q(y) = \sum_{ k \in K} a_k \times \mathbf{Set}(b_k, y) \]
can be written in the form:
 \[ \mathbf{PolyLens}\langle s, t\rangle \langle a, b\rangle = \prod_{k \in K} \mathbf{Set}\left(s_k, \sum_{n \in N} a_n \times \Set(b_n, t_k) \right) \]
Notice that there are two different indexing sets $N$ and $K$. 

There is a coend version of this optic (named, after the compound eyes of insects, ommatidia):
\[ \int^{c_{n k}} 
 \prod_{k \in K} \mathbf{Set} \left(s_k,  \sum_{n \in N} a_n \times c_{n k} \right) \times 
 \prod_{i \in K}  \mathbf{Set} \left(\sum_{m \in N} b_m \times c_{m i}, t_i \right) \]
Here, the coend is taken over the family of sets $c_{nk}$ that are indexed over two different sets, $N$ and $K$, reminiscent of rectangular matrices. 

In order to describe polynomial lens as a compound optic, we'll make a generalization: we'll replace the indexing sets with categories and consider co-presheaves over these categories:
\[a, b \colon [\cat N, \Set] \]
\[ s, t \colon [\cat K, \Set] \]
We interpret the categories $\cat N$ and $\cat K$ as 0-cells in the bicategory $\Cat{Prof}$. We then define a pseudofunctor from $\Cat{Prof}$ to $\Cat{Cat}$. On 0-cells it maps a category $\cat N$ to its co-presheaf category $[\cat N, \Set]$. On 1-cells, it maps a profunctor:
\[ p \colon \cat N^{op} \times \cat K \to \Set \]
to a functor:
\[ c \colon [\cat N, \Set] \to [\cat K, \Set] \]
The action of this functor (here, denoted by $\bullet$) on a co-presheaf $a$ is given by the coend (this is reminiscent of the action of a matrix on a vector):
\[ (c \bullet a) k = \int^{n \colon \cat N} a(n) \times p \langle n, k \rangle \]
Profunctor composition in $\Cat{Prof}$ is defined as:
\[ (p \diamond p') \langle n, k \rangle= \int^m p\langle n, m \rangle \times p'\langle m, k \rangle \]
The action of $p \diamond p'$ on $a$ is isomorphic to the composition of actions:
\[  \int^{n} a(n) \times \int^m p\langle n, m \rangle \times p'\langle m, k \rangle \cong (c' \bullet (c \bullet a)) (k)\]
The unit of profunctor composition is the hom-functor, and the corresponding action is isomorphic,  by co-Yoneda, to identity.

On 2-cells, the pseudofunctor maps natural transformations between profunctors to natural transformations between functors.

In general, compound optics are defined by a pair of pseudofunctors. Here we use the same pseudofunctor twice. The resulting compound optic is given by a coend over 1-cells $p$, which are profunctors $\cat N^{op} \times \cat K \to \Set$:
\[ \int^{p} 
 [\cat K, \Set] \left(s,  \int^{n} a(n) \times p \langle n, - \rangle \right) \times 
 [\cat K, \Set] \left(\int^{n'} b(n') \times p \langle n', - \rangle, t \right) \]
 
 In the special case of $\cat N$ and $\cat K$ being discrete (that is sets of objects), the internal coends turn into coproducts, and the sets of natural transformations into products, resulting in the formula for polynomial optic or ommatidia.
 
 \section{Future Work}
 The current approach can be naturally generalized to the enriched setting.
 
 There is one aspect of simple optics that still resists generalization, namely the profunctor representation. It's not clear how to generalize Tambara modules and Tannakian reconstruction to the new bicategorical setting. 

\bibliographystyle{plain}
\bibliography{CompoundOptics}
\end{document}